\documentclass[11pt]{amsart}
\usepackage{amsbsy,amsfonts}
\usepackage{latexsym,eucal,amsmath,amsthm}
\usepackage{amssymb}
\newcommand{\db}{\mathbb }

\newcommand{\dbR}{{\db R}}

\def\Im{\mbox{\rm Im}}

\def\X{{\bf X}}
\def\R{{\db R}}
\def\cal{\mathcal}

\theoremstyle{plain}
     \newtheorem{theorem}[subsection]{Theorem}
     
     \newtheorem{proposition}[subsection]{Proposition}
     \newtheorem{lemma}[subsection]{Lemma}

\theoremstyle{remark}
     \newtheorem{remark}[subsection]{Remark}

\theoremstyle{definition}
     \newtheorem{definition}[subsection]{Definition}

\begin{document}
\date{}
\title{Global well-posedness for Schr\"odinger equations with
derivative}

\author{J. Colliander}
\thanks{J.E.C. was supported in part by an N.S.F. Postdoctoral Research Fellowship.}
\address{\small University of California, Berkeley}

\author{M. Keel}
\thanks{M.K. was supported in part by N.S.F. Grant DMS 9801558.}
\address{\small Caltech}

\author{G. Staffilani}
\thanks{G.S. was supported in part by N.S.F. Grant DMS 9800879 and the Terman Award.}
\address{\small Stanford University}

\author{H. Takaoka}
\address{\small Tohoku University}

\author{T. Tao}
\thanks{T.T. is a Clay Prize Fellow and was supported in part by grants from the Packard and Sloan Foundations.}
\address{\small University of California, Los Angeles}

\begin{abstract}
We prove that the 1D Schr\"odinger equation with
derivative in the nonlinear term is globally well-posed in $H^{s}$,
for $s>2/3$ for small $L^{2}$ data. The result follows from an
application of the
``I-method''. This method allows  to define a modification of
the  energy norm $H^{1}$ that is  ``almost conserved''
and can be used to perform an iteration argument.
We also remark that
the same argument can be used to prove that
any quintic nonlinear defocusing Schr\"odinger
equation on the line is globally well-posed for large data in  $H^{s}$,
for $s>2/3$ .
\end{abstract}

\maketitle

\section{Introduction}
We consider the derivative nonlinear Schr\"odinger initial value
problem (IVP)
\begin{equation}
\left\{ \begin{array}{l}
i\partial_tu+\partial_{x}^{2}u=i\lambda\partial_{x}(|u|^2u),\\
    u(x,0) = u_{0}(x),\hspace{1.5cm}x \in \dbR, \, t \in \dbR,
\end{array}\right.
\label{ivp}\end{equation}
where $\lambda \in \dbR$. The equation in (1) is a model for the
propagation of circularly polarized Alfv\'en waves in magnetized
plasma with a constant magnetic field \cite{momt, m, ss}.

It is natural to impose the smallness condition
\begin{equation}\label{small-l2}
\| u_0 \|_{L^2} < \sqrt{\frac{2\pi}{\lambda}}
\end{equation}
on the initial data, as this will force the energy to be positive via
the sharp Gagliardo-Nirenberg inequality.  Note that the $L^2$ norm
is conserved by the evolution.

Well-posedness for the  Cauchy problem \eqref{ivp} has been studied
by many authors \cite{hayashi, hayashi-ozawa, hayashi-ozawa2, ozawa,
takaoka:dnls-local, tsutsumi-fukuda, tsutsumi-fukuda-2}. The best 
local well-posedness result is due to
Takaoka~\cite{takaoka:dnls-local}, where a gauge transformation and
the Fourier restriction method is used to obtain local well-posedness
in  $H^{s}, \, \, s\geq 1/2$.  In \cite{takaoka:dnls-global} Takaoka
showed this result is sharp in the sense that the data map fails to
be $C^{3}$ or uniformly $C^0$ for $s < 1/2$ (cf. Bourgain
\cite{borg:measures} and Biagioni-Linares \cite{bl}).

In \cite{ozawa} global well-posedness is obtained for \eqref{ivp} in
$H^1$ assuming the smallness condition \eqref{small-l2}. The argument
there is based on two gauge
transformations performed in order to remove the derivative in the
nonlinear term.  This was improved by
Takaoka~\cite{takaoka:dnls-global}, who proved global well-posed in
$H^{s}$ for $s>\frac{32}{33}$ assuming \eqref{small-l2}. The method
of proof is based on the idea of Bourgain~\cite{borg:refinements,
borg:book} of estimating separately the evolution of low frequencies
and of high frequencies of the initial data.

In this paper we improve the global well-posedness result further:

\begin{theorem}\label{main}  The Cauchy problem \eqref{ivp} is
globally well-posed in $H^s$ for $s > 2/3$, assuming the smallness
condition \eqref{small-l2}.
\end{theorem}

The proof of Theorem \ref{main} is based on the ``I-method'' used by
the authors in other non-linear Cauchy problems in \cite{keel:mkg},
\cite{ckstt:1}, \cite{ckstt:2}, \cite{ckstt:3} (see also
\cite{keel:wavemap}).  The basic idea is as follows.  After a
rescaling, we define a new energy $E_N(u)(t)$ for the solution $u$
that depends on a parameter $N \gg 1$. We prove a local well-posed
result in the norm associated to $E_N$ on intervals of length $\sim
1$, and finally we perform an iteration on  the time intervals. The
reason why this iteration can be globally extended is
that the increment of the energy $E_{N}(u)(t)$
over each time interval is very small. In other words the argument
is successful because the energy  $E_{N}(u)(t)$ is {\em almost
conserved}.

After the proof of Theorem \ref{main} is completed, we will
briefly remark that using the same techniques one can also show
that the
1D defocusing quintic nonlinear Schr\"odinger is global well-posed
for initial data in $H^{s}, s>2/3$. The details of the proof of this
fact will appear in a different paper.

The restriction $s > 2/3$ is probably not sharp, and might be
improvable either by more sophisticated multilinear estimates and
better estimates on the symbols $M_4$, $M_6$, $M_8$ which appear in
our argument, or by using the ``correction term'' strategy of
\cite{ckstt:2}.  In fact one may reasonably conjecture that one could
extend the global well-posedness result to match the local result at
$s > 1/2$.  We will not pursue these matters here.

\section{Notation}

To prove Theorem \ref{main} we may assume $2/3 < s < 1$, since for
the $s \geq  1$ the result is contained in \cite{ozawa,
takaoka:dnls-global}.  Henceforth $2/3 < s
< 1$ shall be fixed.  Also, by rescaling $u$, 
we may assume $\lambda = 1$.

We use $C$ to denote various constants depending on $s$; if $C$
depends on other quantities as well, this will be indicated by
explicit subscripting, e.g. $C_{\|u_0\|_2}$ will depend on both $s$
and $\|u_0\|_2$.  We use $A \lesssim B$ to denote an estimate of the
form $A \leq C B$.  We use $a+$ and $a-$ to denote expressions of the
form $a+\varepsilon$ and $a-\varepsilon$, where $0 < \varepsilon \ll
1$ depends only on $s$.

We use $\|f\|_{p}$ to denote the $L^{p}(\dbR)$ norm, and $L^q_t
L^r_x$ to denote the mixed norm
$$ \| f\|_{L^q_t L^r_x} := (\int \|f(t)\|_r^q\ dt)^{1/q}$$
with the usual modifications when $q=\infty$.

We define the spatial Fourier transform of $f(x)$ by
$$ \hat f(\xi) := \int_{\R} e^{-i x \xi} f(x)\ dx$$
and the spacetime Fourier transform $u(t,x)$ by
$$ \tilde u(\tau, \xi) := \int_{\R}\int_{\R} e^{-i (x \xi + t \tau)}
u(t,x)\ dt dx.$$
Note that the derivative $\partial_x$ is conjugated to multiplication
by $i\xi$ by the Fourier transform.

We shall also define $D_x$ to be the Fourier multiplier with symbol
$\langle \xi \rangle := 1 + |\xi|$.  We can then define the Sobolev
norms $H^s$ by
$$ \| f \|_{H^s} := \| D_x^s f\|_2 = \| \langle \xi \rangle^s \hat f
\|_{L^2_\xi}.$$
We also define the spaces $X^{s,b}(\R \times \R)$ (first introduced
in \cite{borg:xsb}) on $\R \times \R$ by
$$ \| u \|_{X^{s,b}(\R \times \R)} := \| \langle \xi \rangle^s
\langle \tau - |\xi|^2 \rangle^b \hat{u}(\xi,\tau)\|_{L^2_\tau L^2_\xi}.$$
We often abbreviate $\| u\|_{s,b}$ for $\| u \|_{X^{s,b}(\R \times
\R)}$.
For any time interval $I$, we define the restricted spaces $X^{s,b}(I
\times \R)$ by
$$ \| u \|_{X^{s,b}(I \times \R)} := \inf \{ \| U \|_{s,b}:
U|_{I \times \R} = u \}.$$
We shall take advantage of the Strichartz estimates
\begin{equation}\label{strichartz-6}
\| u \|_{L^6_t L^6_x} \lesssim \| u \|_{0,1/2+}
\end{equation}
and
\begin{equation}\label{strichartz-2}
\| u \|_{L^\infty_t L^2_x} \lesssim \| u \|_{0,1/2+}
\end{equation}
(see e.g. \cite{borg:xsb}).  From \eqref{strichartz-2} and Sobolev embedding we
observe
\begin{equation}\label{strichartz-infty}
\| u \|_{L^\infty_t L^\infty_x} \lesssim \| u \|_{1/2+,1/2+}
\end{equation}
In our arguments we shall be using the trivial embedding
$$ \| u \|_{s_1,b_1} \lesssim \| u \|_{s_2,b_2} \hbox{ whenever } s_1
\leq s_2, b_1 \leq b_2$$
so frequently that we will not mention this embedding explicitly.

We now give some useful notation for multilinear expressions.  If $n
\geq 2$ is an even integer, we define a \emph{(spatial) multiplier of
order $n$} to be any function $M_n(\xi_1, \ldots, \xi_n)$ on the
hyperplane
$$ \Gamma_n := \{ (\xi_1, \ldots, \xi_n) \in \R^n: \xi_1 + \ldots +
\xi_n = 0 \},$$
which we endow with the standard measure $\delta(\xi_1 + \ldots +
\xi_n)$, where $\delta$ is the Dirac delta.

If $M_n$ is a multiplier of order $n$ and $f_1, \ldots, f_n$ are
functions on $\R$, we define the quantity $\Lambda_n(M_n;f_1, \ldots,
f_n)$ by
$$ \Lambda_n(M_n;f_1, \ldots, f_n)
:= \int_{\Gamma_n} M_n(\xi_1, \ldots, \xi_n) \prod_{j=1}^n \hat
f_j(\xi_j).$$
We adopt the notation
$$ \Lambda_n(M_n;f) := \Lambda_n(M_n; f, \bar{f}, f, \bar{f}, \ldots,
f, \bar{f}).$$
Observe that $\Lambda_n(M_n;f)$ is invariant under permutations of
the even $\xi_j$ indices, or of the odd $\xi_j$ indices.

If $M_n$ is a multiplier of order $n$, $1 \leq j \leq n$ is an index,
and $k \geq 1$ is an even integer, we define the {\em elongation}
$\X^k_j(M_n)$ of $M_n$ to be the multiplier of order $n+k$ given by
$$ \X^k_j(M_n)(\xi_1, \ldots, \xi_{n+k})
:= M_n(\xi_1, \ldots, \xi_{j-1}, \xi_j + \ldots + \xi_{j+k},
\xi_{j+k+1}, \ldots, \xi_{n+k}).$$
In other words, $\X^k_j$ is the multiplier obtained by replacing
$\xi_j$ by $\xi_j + \ldots + \xi_{j+k}$ and advancing all the indices
after $\xi_j$ accordingly.

We shall often write $\xi_{ij}$ for $\xi_i + \xi_j$, $\xi_{ijk}$ for
$\xi_i + \xi_j + \xi_k$, etc.  We also write $\xi_{i-j}$ for $\xi_i -
\xi_j$, $\xi_{ij-klm}$ for $\xi_{ij} - \xi_{klm}$, etc.

\section{The Gauge Transformation and the conservation laws}

In this section we apply the gauge transform used in \cite{ozawa} in
order to improve the derivative nonlinearity.

\begin{definition}\label{gauge-def}  We define the non-linear map
${\cal G}: L^2(\R) \to L^2(\R)$ by
$$
{\cal G} f(x) := e^{-i\int_{-\infty}^{x}|f(y)|^{2}dy}f(x).$$
The inverse transform ${\cal G}^{-1} f$ is then given by
$$
{\cal G}^{-1} f(x) := e^{i\int_{-\infty}^{x}|f(y)|^{2}dy}f(x).$$
\end{definition}

This transform is well behaved on $H^s$:

\begin{lemma}\label{bicontinuity}
The map ${\cal G}$ is a bicontinuous map from $H^s$ to $H^s$.
\end{lemma}

A similar statement holds for $0 \leq s \leq 1/2$, but we shall not
need it here.

\begin{proof}
We shall just prove the continuity of ${\cal G}$, as the continuity
of ${\cal G}^{-1}$ is proven similarly.

Define $Lip$ to be the space of functions with norm
$$ \| f\|_{Lip} := \| f\|_{\infty} + \| f' \|_{L^\infty}.$$
Since $s>1/2$, we see from Sobolev embedding that the nonlinear map
$f \mapsto e^{-i\int_{-\infty}^{x}|f(y)|^{2}dy}$ continuously maps
$H^s$ to $Lip$.  It therefore suffices to show the product estimate
$$ \| fg \|_{H^s} \lesssim \| f\|_{H^s} \| g \|_{Lip}.$$
But this estimate follows immediately from the Leibnitz rule and
H\"older when $s=0$ or $s=1$, and the intermediate cases then follow
by interpolation. \end{proof}

Set $w_0 := {\cal G} u_0$, and $w(t) := {\cal G} u(t)$ for all times
$t$.  A straightforward calculation shows that the IVP \eqref{ivp}
can be transformed to
\begin{equation}
\left\{ \begin{array}{l}
i\partial_tw+\partial_{x}^{2}w=-i w^{2}\partial_{x}\bar{w}
-\frac{1}{2}|w|^{4}w,\\
    w(x,0) = w_{0}(x),\hspace{1.5cm}x \in \dbR, \, t \in \dbR,
\end{array}\right.
\label{givp1}\end{equation}
Also, the smallness condition \eqref{small-l2} becomes
\begin{equation}\label{small-l2-w}
\| w_0 \|_{L^2} < \sqrt{2\pi}.
\end{equation}
By Lemma \ref{bicontinuity} we thus see that global well-posedness of
\eqref{ivp} in $H^s$ is equivalent to that of \eqref{givp1}.  From
\cite{ozawa, takaoka:dnls-local, takaoka:dnls-global},  we know that 
both Cauchy problems are
locally well-posed in $H^s$ and globally well-posed in $H^1$ assuming
\eqref{small-l2-w}.  By standard limiting arguments, we thus see that
Theorem \ref{main} will follow if we can show

\begin{proposition}\label{global}  Let $w$ be a global $H^1$ solution
to \eqref{givp1} obeying \eqref{small-l2-w}.  Then for any $T > 0$ we
have
$$ \sup_{0 \leq t \leq T} \| w(t) \|_{H^s} \lesssim
C_{\|w_0\|_{H^s},\|w_0\|_2,T}$$
where the right-hand side does not depend on the $H^1$ norm of $w$.
\end{proposition}

Just by looking at the equation in (\ref{givp1}) it is not easy to
understand why
this should be  better than the equation in (\ref{ivp}). In fact we
still
  see a derivative, and moreover a quintic nonlinearity has been
introduced.
But it was made clear in
\cite{ozawa, hayashi-ozawa3, takaoka:dnls-local} how a derivative of 
the complex
conjugate of the
solution $w$ can be handled while a derivative of $w$ cannot.
Also the
quintic term is not  going to introduce any extra trouble.

Let $n \geq 2$ be an even integer, and let $M_n$ be a multiplier of
order $n$.
 From \eqref{givp1} we have
$$ \partial_t w =  i w_{xx} - w \bar w_x w + \frac{i}{2} w \bar{w} w
\bar{w} w$$
and
$$ \partial_t \bar{w} = -i w_{xx} - \bar{w} w_x \bar{w} - \frac{i}{2}
\bar{w} w \bar{w} w \bar{w}.$$
Taking the Fourier transform of these identities, we obtain the
useful differentiation law
\begin{equation}\label{diff}
\begin{split}
\partial_t \Lambda_n(M_n;w(t))
&=
i \Lambda_n(M_n \sum_{j=1}^n (-1)^{j} \xi_j^2; w(t))\\
&- i \Lambda_{n+2}(\sum_{j=1}^n \X^2_j(M_n) \xi_{j+1}; w(t))\\
&+\frac{i}{2} \Lambda_{n+4}(\sum_{j=1}^n (-1)^{j-1} \X^4_j(M_n); w(t))
\end{split}
\end{equation}
for any even integer $n \geq 2$ and any multiplier $M_n$ of order $n$.

We now turn to the conservation laws that the solution $w$ of
(\ref{givp1}) enjoys.  What follows  in this section was originally
described by
Ozawa in \cite{ozawa}, however we have redone the computations in our
own notation as this will prove useful later.

\begin{definition}\label{energy-def}  If $f \in H^1(\R$), we define
the energy $E(f)$ by
$$ E(f) := \int \partial_x f\partial_{x}\overline{f} \ dx
- \frac{1}{2}\Im
\int f \overline{f} f \partial_x \overline{f}\ dx.$$
\end{definition}

By Plancherel, we may write $E(f)$ using the $\Lambda$ notation as
$$ E(f) = -\Lambda_2(\xi_1\xi_2;f) -  \frac{1}{2} \Im
\Lambda_4(i\xi_4;f).$$
Expanding out the second term using $\Im(z) = (z - \bar{z})/2i$, and
using symmetry, we may rewrite this as
\begin{equation}\label{energy-lambda}
E(f) = -\Lambda_2(\xi_1\xi_2;f) + \frac{1}{8} \Lambda_4(\xi_{13-24};f).
\end{equation}

\begin{lemma}\label{conservation}\cite{ozawa}  If $w$ is an $H^1$
solution to \eqref{givp1} for times $t \in [0,T]$, then we have
$$ \| w(t) \|_2 = \| w_0 \|_2$$
and
$$ E(w(t)) = E(w_0)$$
for all $t \in [0,T]$.
\end{lemma}

\begin{proof}
These conservation laws are proven in \cite{ozawa}, however we give a
proof based on the identity \eqref{diff}, as the proof here will be
needed later on.

We of course have
$$ \| w(t)\|_2^2 = \Lambda_2(1; w(t)).$$
In the rest of this proof we shall drop the $w(t)$ from the $\Lambda$
notation.
Differentiating the previous and applying \eqref{diff}, we obtain
$$ \partial_t \| w(t) \|_2^2
=
-i \Lambda_2(\xi_1^2 - \xi_2^2) - i \Lambda_4(\xi_2 + \xi_3)
+ \frac{i}{2} \Lambda_6(1-1+1-1+1-1).$$
The first term vanishes since $\xi_{12} = 0$.  The second term can be
symmetrized to $-\frac{i}{2} \Lambda_4(\xi_{1234})$ which vanishes.
The third term clearly vanishes.  This proves the $L^2$ conservation.

Now we prove energy conservation.  From (\ref{energy-lambda})
we have
\begin{equation}
\partial_{t}E(t)=-\partial_{t}\Lambda_{2}(\xi_{1}\xi_{2})
+\frac{1}{8}\partial_{t}\Lambda_{4}(\xi_{13-24})
\label{dte}\end{equation}
and from \eqref{diff}
$$
\partial_t \Lambda_2(\xi_1\xi_2) =
-i \Lambda_2(\xi_1 \xi_2(\xi_1^2 - \xi_2^2))
- i \Lambda_4(\xi_{123} \xi_4 \xi_2 + \xi_1 \xi_{234} \xi_3)
+ \frac{i}{2} \Lambda_6(\xi_{12345} \xi_6 - \xi_1 \xi_{23456}).$$
The $\Lambda_2$ term vanishes since $\xi_{12} = 0$.  To simplify the
$\Lambda_4$ term we write $\xi_{123} = - \xi_4$, $\xi_{234} = -\xi_1$
and then symmetrize.  To simplify the $\Lambda_6$ term we write
$\xi_{12345} = - \xi_6$, $\xi_{23456} = -\xi_1$ and then symmetrize,
to obtain
$$
\partial_t \Lambda_2(\xi_1\xi_2) =
\frac{i}{2}
\Lambda_4(\xi_1^2 \xi_3 + \xi_2^2 \xi_4 + \xi_3^2 \xi_1 + \xi_4^2
\xi_2)
+
\frac{i}{6} \Lambda_6(\xi_1^2 - \xi_2^2 + \xi_3^2 - \xi_4^2 + \xi_5^2
- \xi_6^2).
$$
We may simplify the $\Lambda_4$ term further, using the identity
\begin{align*}
\xi_1^2 \xi_3 + \xi_2^2 \xi_4 + \xi_3^2 \xi_1 + \xi_4^2 \xi_2
&= \xi_1 \xi_3 \xi_{13} + \xi_2 \xi_4 \xi_{24}\\
&= \xi_{13} (\xi_1 \xi_3 - \xi_2 \xi_4)\\
&= \xi_{13} (-\xi_1 \xi_{124} - \xi_2 \xi_4)\\
&= -\xi_{13} (\xi_1 + \xi_2) (\xi_1 + \xi_4)\\
&= -\xi_{12} \xi_{13} \xi_{14}
\end{align*}
to obtain
\begin{equation}\label{l4-d}
\partial_t \Lambda_2(\xi_1\xi_2) =
-\frac{i}{2}
\Lambda_4(\xi_{12} \xi_{13} \xi_{14})
+ \frac{i}{6} \Lambda_6(\xi_1^2 - \xi_2^2 + \xi_3^2 - \xi_4^2 +
\xi_5^2 - \xi_6^2).
\end{equation}

We now consider the second component of the energy.  From
\eqref{diff} we have
\begin{align*}
\partial_t \Lambda_4(\xi_{13-24}) &=
i \Lambda_4(\xi_{13-24} (\xi_1^2 - \xi_2^2 + \xi_3^2 - \xi_4^2))\\
&- i \Lambda_6(
\xi_{1235-46}\xi_2
+ \xi_{15-2346} \xi_3
+ \xi_{1345-26}\xi_4
+ \xi_{13-2456}\xi_5)\\
&+ \frac{i}{2} \Lambda_8(\xi_{123457-68} - \xi_{17-234568} +
\xi_{134567-28} - \xi_{13-245678}).
\end{align*}
The $\Lambda_8$ term symmetrizes to $i \Lambda_8(\xi_{12345678})$
which vanishes.
The $\Lambda_6$ term can be rewritten as
$$ 2i \Lambda_6(\xi_{46} \xi_2 - \xi_{15}\xi_3 + \xi_{26} \xi_4 -
\xi_{13}\xi_5)$$
which we rewrite as
$$ 2i \Lambda_6(\xi_{246} \xi_2 - \xi_{135}\xi_3 + \xi_{246} \xi_4 -
\xi_{135}\xi_5)
- 2i \Lambda_6(\xi_2^2 - \xi_3^2 + \xi_4^2 - \xi_5^2).$$
The first term symmetrizes to $\frac{4i}{3}\Lambda_{6}(\xi_{246}^2 -
\xi_{135}^2)$
which vanishes.  The second term symmetrizes to
$$ \frac{4i}{3} \Lambda_6(\xi_1^2 - \xi_2^2 + \xi_3^2 - \xi_4^2 +
\xi_5^2 - \xi_6^2).$$
Finally, consider the $\Lambda_4$ term.  We may factorize
$$\xi_{13-24} (\xi_1^2 - \xi_2^2 + \xi_3^2 - \xi_4^2)
= \xi_{13-24} (\xi_{1-2} \xi_{12} + \xi_{3-4} \xi_{34}).$$
Since $\xi_{12} = \xi_{-34}$ and $\xi_{13} = -\xi_{24}$, we may
simplify this as
$$ 2\xi_{13} \xi_{12} (\xi_{1-2} - \xi_{3-4})
= 4 \xi_{12} \xi_{13} \xi_{14}.$$
Combining all these identities we thus have
$$
\frac{1}{8}\partial_t \Lambda_4(\xi_{13-24}) =
-\frac{1}{2}i \Lambda_4(\xi_{12} \xi_{13} \xi_{14}) - \frac{i}{6}
\Lambda_6(\xi_1^2 - \xi_2^2 + \xi_3^2 - \xi_4^2 + \xi_5^2 -
\xi_6^2).$$
Combining this with \eqref{l4-d} and (\ref{dte}) we obtain
$$ \partial_t E(w(t)) = 0$$
and the claim follows.
\end{proof}

Heuristically, the energy $E(w(t))$ has the same strength as $\|
w(t)\|_{H^1}^2$.  We can make this precise:

\begin{lemma}\label{gagliardo}  Let $f$ be an $H^1$ function on $\R$
such that
$\|f\|_{2}<\sqrt{2\pi}$.   Then we have
\begin{equation}\label{gee}
\|\partial_{x}f\|_2\leq C_{\|f\|_2} E(f)^{1/2}
\end{equation}
where $C_{\|f\|_2}$ depends only on $\|f\|_2$.
\label{ee}\end{lemma}
\begin{proof}
Define the function
$$ g(x) := \exp(i\frac{3}{4} \int_{-\infty}^x |f(y)|^2\ dy) f(x).$$
A routine computation shows that
$$ \|g\|_2 = \|f\|_2 < \sqrt{2\pi}$$
and
$$ E(f) = \| \partial_x g \|_2^2 - \frac{1}{16} \|g\|_6^6.$$
 From the sharp Gagliardo-Nirenberg inequality \cite{w}
\begin{equation}\label{gn}
\|g\|_{6}^{6}\leq \frac{4}{\pi^{2}}\|g\|_{2}^{4}
\|\partial_{x}g\|_{2}^{2}
\end{equation}
we therefore have
$$ \| \partial_x g \|_2 \lesssim C_{\|f\|_2} E(f)^{1/2}.$$
 From the definition of $g$ we have
$$ f(x) = \exp(-i\frac{3}{4} \int_{-\infty}^x |g(y)|^2\ dy) g(x)$$
and so we have
$$ \| \partial_x f\|_2 \lesssim \|\partial_x g\|_2 + \| g^3\|_2.$$
By another application of \eqref{gn} we thus obtain \eqref{gee}.
\end{proof}

\section{The Almost Conserved Energy Norm}
It remains to prove Proposition \ref{global}.  Fix $w$, $T$.  We also
let $N \gg 1$ be a large parameter depending on $T$, $\|w_0\|_2$, and
$\|w_0\|_{H^s}$ which we shall choose later.

Because we do not want to use the $H^1$ norm of $w$, we cannot
directly use the energy $E(w(t))$ defined
above. So we are looking for a substitute notion of ``energy''
that can be
defined for a less regular solution and that has a very slow increment
in time. In the frequency space let us consider an even  $C^{\infty}$
monotone multiplier
$m(\xi)$ taking values in $[0,1]$ such that
\begin{equation}
m(\xi):=\left\{ \begin{array}{l}
1, \, \, \, \mbox{ if } \, |\xi|<N,\\
\left(\frac{|\xi|}{N}\right)^{s-1}\, \, \, \mbox{ if } \, |\xi|> 2N.
\end{array}\right.
\label{mul}\end{equation}
We define the multiplier operator $I:H^{s}\longrightarrow H^{1}$
such that $\widehat{Iw}(\xi):=m(\xi)\widehat{w}(\xi)$.  This operator
is smoothing of order $1-s$; indeed we have
\begin{equation}\label{i-smoothing}
\| u\|_{s_0,b_0} \lesssim \| Iu \|_{s_0 + 1-s, b_0} \lesssim N^{1-s}
\| u \|_{s_0,b_0}
\end{equation}
for any $s_0, b_0 \in \R$.

Our substitute energy will be defined by
$$ E_N(w) := E(Iw).$$
Note that this energy makes sense even if $w$ is only in $H^s$.

In general the energy $E_N(w(t))$ is not conserved in time, but we
will show that the increment is very small in terms of $N$.  This
will be accomplished in three stages.  First, in Proposition
\ref{energy-increment} below, we write the increment of $E_N(w(t))$
as a multilinear expression in $w$.  Then, in Lemma \ref{pinc}, we
estimate these multilinear expressions in terms of the norm $\| I
w\|_{1,1/2+}$, gaining a power of $N^{-1+}$ in the process.  Finally,
in Theorem \ref{local} (and Lemma \ref{gagliardo}), we control the
norm $\| Iw\|_{1,1/2+}$ back in terms of $E_N(w(t))$.

\begin{proposition}\label{energy-increment}
Let $w$ be an $H^1$ global solution to \eqref{givp1}.  Then for any
$T \in \R$ and $\delta > 0$ we have
$$
E_N(w(T+\delta)) - E_N(w(T)) =
\int_T^{T+\delta} [\Lambda_4(M_4;w(t)) + \Lambda_6(M_6;w(t)) +
\Lambda_8(M_8;w(t))]\ dt
$$
where the multipliers $M_4$, $M_6$, $M_8$ are given by
\begin{align*}
M_4
&:=
C_1 m_1 m_2 m_3 m_4 \xi_{12} \xi_{13} \xi_{14} + C_2
(m_1^2 \xi_1^2 \xi_3 + m_2^2 \xi_2^2 \xi_4 + m_3^2 \xi_3^2 \xi_1 +
m_4^2 \xi_4^2 \xi_2) \\
M_6 &:=  C_3 \sum_{j=1}^6 (-1)^{j-1} m_j^2 \xi_j^2
+ C_4 \sum_{\{a,c,e\}=\{1,3,5\}, \{b,d,f\}=\{2,4,6\}}
m_a m_b m_c m_{def} \xi_{ac} \xi_e - m_{abc} m_d m_e m_f \xi_{df}
\xi_b
\\
M_8 &:= C_5 \sum_{\{a,c,e,g\}=\{1,3,5,7\}; \{b,d,f,h\} = \{2,4,6,8\}}
m_a m_b m_c m_{defgh} \xi_{ac-bdefgh} - m_{abcde} m_f m_g m_h
\xi_{abcdeg-fh}
\end{align*}
where $C_1, \ldots, C_5$ are absolute constants and we adopt the
abbreviations $m_i$ for $m(\xi_i)$, $m_{ij}$ for $m(\xi_{ij})$, etc.
Furthermore, if $|\xi_j| \ll N$ for all $j$, then the multipliers
$M_4$, $M_6$, $M_8$ all vanish.
\end{proposition}

\begin{proof}
 From \eqref{energy-lambda} we have
$$ E_N(w(t)) =
-\Lambda_2(m_1 \xi_1 m_2 \xi_2;w(t)) + \frac{1}{8}
\Lambda_4(\xi_{13-24} m_1 m_2 m_3 m_4;w(t)).$$
Henceforth we omit the $w(t)$ from the $\Lambda$ notation.  By
\eqref{diff} we have
\begin{align*}
\partial_t \Lambda_2(m_1 \xi_1 m_2 \xi_2)
&= -i \Lambda_2(m_1 \xi_1 m_2 \xi_2 (\xi_1^2 - \xi_2^2))\\
&- i \Lambda_4(m_{123} \xi_{123} m_4 \xi_4 \xi_2 + m_1 \xi_1 m_{234}
\xi_{234} \xi_3)\\
&+ \frac{i}{2} \Lambda_6(m_{12345}\xi_{12345} m_6 \xi_6
- m_1 \xi_1 m_{23456} \xi_{23456}).
\end{align*}
The $\Lambda_2$ term vanishes since $\xi_{12} = 0$.  To simplify the
$\Lambda_4$ term, we use $\xi_{123} = -\xi_4$ and $\xi_{234} =
-\xi_1$ and then symmetrize to obtain the second term of $M_4$.  To
simplify the $\Lambda_6$ term, we use $\xi_{12345} = -\xi_6$ and
$\xi_{23456} = - \xi_1$ and then symmetrize to get the first term of
$M_6$.

In a similar vein we have
\begin{align*}
\partial_t \Lambda_4(\xi_{13-24} m_1 m_2 m_3 m_4)
&= -i \Lambda_4(\xi_{13-24} m_1 m_2 m_3 m_4 (\xi_1^2 - \xi_2^2 +
\xi_3^2-\xi_4^2) )\\
&- i \Lambda_6(\xi_{1235-46} m_{123} m_4 m_5 m_6 \xi_2
+ \xi_{15-2346} m_1 m_{234} m_5 m_6 \xi_3\\
& \quad + \xi_{1345-26} m_1 m_2 m_{345} m_6 \xi_4
+ \xi_{13-2456} m_1 m_2 m_3 m_{456} \xi_5)\\
&+ \frac{i}{2} \Lambda_8(\xi_{123457-68} m_{12345} m_6 m_7 m_8 -
\xi_{17-234568} m_1 m_{23456} m_7 m_8 \\
& \quad + \xi_{134567-28} m_1 m_2 m_{34567} m_8 -
\xi_{13-245678} m_1 m_2 m_3 m_{45678}).
\end{align*}
The $\Lambda_4$ term is of the form of the first term of $M_4$, by
the argument used to prove \eqref{l4-d}.  To simplify the $\Lambda_6$
term, we use $\xi_{1235-46} = -2\xi_{46}$ and similarly for the other
four terms, then symmetrize to obtain the second term of $M_6$.
Finally if we symmetrize the $\Lambda_8$ term we obtain $M_8$.  The
first part of the Proposition then follows from the Fundamental
Theorem of Calculus applied to the function
$t\longrightarrow E_{N}(w(t))$.

If all the frequencies are $\ll N$, then all the $m_i$, $m_{ij}$,
etc. terms are equal to 1.  In this case our calculations are
identical to those in Lemma \ref{conservation} and so our symbols
$M_4$, $M_6$, $M_8$ will vanish by the computations given in that
Lemma.
\end{proof}

\section{Local estimates}
In Lemma \ref{pinc} we shall estimate the expression in Proposition
\ref{energy-increment}.  It turns out that one cannot estimate this
expression effectively just by using spatial norms such as $\|
Iw\|_{H^1}$ (as is done for some simple equations in
\cite{borg:book}), but one must use spacetime norms such as $\|
Iw\|_{1,1/2+}$.  The purpose of this section is to obtain the
required control on these spacetime norms:

\begin{theorem}\label{local}  Let $w$ be a $H^1$ global solution to
\eqref{givp1}, and let $T \in \R$ be such that
$$ \| Iw(T) \|_{H^1} \leq C_0$$
for some $C_0 > 0$.  Then we have
$$ \| Iw \|_{X^{1,1/2+}([T,T + \delta] \times \R)} \lesssim 1$$
for some $\delta > 0$ depending on $C_0$.
\label{lwp}\end{theorem}
We now prove Theorem \ref{lwp}.  We shall be able to exploit the
estimates in \cite{takaoka:dnls-local}.  By standard iteration
arguments (see e.g. \cite{borg:xsb}, \cite{kpv:xsb}, \cite{kpv:kdv},
\cite{takaoka:dnls-local}, \cite{takaoka:dnls-global}) it suffices to
prove
\begin{lemma}
We have

\begin{align}
   \|I(w_{1}\partial_{x}\overline{w_{2}} w_{3})\|_{X_{1,b-1}(\R \times
\R)}&\lesssim
   \prod_{i=1}^{3}\|Iw_{i}\|_{X_{1,1/2+}(\R \times
\R)}\label{4-linear}\\
   \|I(w_{1}\overline{w_{2}}w_{3}\overline{w_{4}}w_{5})\|_{X_{1,b-1}
   (\R \times \R)}
   &\lesssim
   \prod_{i=1}^{5}\|Iw_{i}\|_{X_{1,1/2+}(\R \times
\R)}.\label{6-linear}
\end{align}
for all Schwarz functions $w_i$ and some $b > 1/2$ (in fact we may
take any $1/2 < b < 5/8$).
\end{lemma}
\begin{proof}
By Plancherel and duality it suffices to show
$$
|\int_* \frac{ m(\xi_4) \langle \xi_4 \rangle \langle \tau_{4}+ \xi_4^2
\rangle^{b-1} \xi_2 }
{ \prod_{j=1}^3 m(\xi_j) \langle \xi_j \rangle \langle \tau_j -
(-1)^{j-1}\xi_j^2 \rangle^{1/2+} } \prod_{j=1}^4 F_j(\tau_j,\xi_j)|
\lesssim \prod_{j=1}^4 \|F_j\|_{L^2_{\tau_j} L^2_{\xi_j}}$$
and
$$|\int_{**} \frac{ m(\xi_6) \langle \xi_6 \rangle \langle \tau_{6} +
\xi_6^2 \rangle^{b-1} }
{ \prod_{j=1}^5 m(\xi_j) \langle \xi_j \rangle \langle \tau_j -
(-1)^{j-1}\xi_j^2 \rangle^{1/2+} }
\prod_{j=1}^6 F_j(\tau_j,\xi_j)|
\lesssim \prod_{j=1}^6 \|F_j\|_{L^2_{\tau_j} L^2_{\xi_j}}$$
for all functions $F_1, \ldots, F_6$, where $\int_*$, $\int_{**}$
denotes integration over the measure $\delta(\tau_1 + \ldots +
\tau_4) \delta(\xi_1 + \ldots + \xi_4)$ and $\delta(\tau_1 + \ldots +
\tau_6) \delta(\xi_1 + \ldots + \xi_6)$ respectively.

We may assume that the $F_j$ are all real and non-negative.  We now
observe the pointwise estimate
$$ \frac{ m(\xi_n) \langle \xi_n \rangle^{1-s} }
{ \prod_{j=1}^{n-1} m(\xi_j) \langle \xi_j \rangle^{1-s} } \lesssim
1$$
for $n=4,6$ and all $\xi_1, \ldots, \xi_n$ such that $\xi_1 + \ldots
+ \xi_n = 0$.  To see this, we use symmetry to assume that $|\xi_1|
\geq \ldots \geq |\xi_{n-1}|$, so that $|\xi_n| \lesssim |\xi_1|$.
Since $m(\xi) \langle \xi \rangle^{1-s}$ is essentially increasing in
$|\xi|$, we thus see that
$$ \frac{ m(\xi_n) \langle \xi_n \rangle^{1-s} }
{ \prod_{j=1}^{n-1} m(\xi_j) \langle \xi_j \rangle^{1-s} } \lesssim
\frac{ 1}
{ \prod_{j=2}^{n-1} m(\xi_j) \langle \xi_j \rangle^{1-s} }.$$
Since $m(\xi) \langle \xi \rangle^{1-s} \gtrsim 1$ for all $\xi$, the
claim follows.

Because of this estimate, we only need to show the estimates
\begin{equation}\label{est-4}
|\int_*  \frac{ \langle \xi_4 \rangle^s \langle \tau_{4} + \xi_4^2
\rangle^{b-1} \xi_2 }
{ \prod_{j=1}^3 \langle \xi_j \rangle^s \langle \tau_j -
(-1)^{j-1}\xi_j^2 \rangle^{1/2+} } \prod_{j=1}^4 F_j(\tau_j,\xi_j)|
\lesssim \prod_{j=1}^4 \|F_j\|_{L^2_{\tau_j} L^2_{\xi_j}}
\end{equation}
and
\begin{equation}\label{est-6}
|\int_{**}
\frac{ \langle \xi_6 \rangle^s \langle \tau_{6} + \xi_6^2 \rangle^{b-1} }
{ \prod_{j=1}^5 \langle \xi_j \rangle^s \langle \tau_j -
(-1)^{j-1}\xi_j^2 \rangle^{1/2+} } \prod_{j=1}^6 F_j(\tau_j,\xi_j)|
\lesssim \prod_{j=1}^6 \|F_j\|_{L^2_{\tau_j} L^2_{\xi_j}}.
\end{equation}
The estimate \eqref{est-4} is equivalent to the first estimate of
Lemma 3.1 in \cite{takaoka:dnls-global} after undoing the duality and
Plancherel, so it suffices to prove \eqref{est-6}.  By undoing the
duality we can write this as
$$ \| w_1 \overline{w_2} w_3 \overline{w_4} w_5 \|_{s,b-1}
\lesssim \prod_{j=1}^5 \| w_j \|_{s,1/2+}.$$
We may assume that the Fourier transforms $\tilde w_j$ are all real
and non-negative.
By using $|\xi_{1}+\xi_{2}+\xi_{3}+\xi_{4}+\xi_{5}|^{s}\lesssim
\sum_{i=1}^{5}|\xi_{i}|^{s}$, it suffices to prove estimates of the
form
$$ \| (D_x^s w_1) \overline{w_2} w_3 \overline{w_4} w_5 \|_{0,b-1}
\lesssim \prod_{j=1}^5 \| w_j \|_{s,1/2+},$$
plus similar estimates when $D_x^s$ falls on one of the other
functions.  We shall only prove the displayed estimate, as the others
are similar.  We may estimate the $X^{0,b-1}$ norm by the $L^2_t
L^2_x$ norm.  But then the claim follows from three applications of
\eqref{strichartz-6} and two applications of
\eqref{strichartz-infty}, and H\"older (ensuring that the term with
the $D_x^s$ is estimated using \eqref{strichartz-6}).
\end{proof}

\section{Proof of Proposition \ref{global}}\label{global-sec}
We can now prove Proposition \ref{global}, which as remarked before
will give Theorem \ref{main}. Let $T$, $w$ be as in the Proposition.
Our constants may depend on $\|w_0\|_2$ and $\|w_0\|_{H^s}$.

We start by rescaling the solution $w$.  Let $\mu > 0$ be chosen
later.  We observe that
$w$ is a solution for the IVP (\ref{givp1}) if and only if
\[w^{\mu}(t,x)=\frac{1}{\mu^{1/2}}
w\left(\frac{t}{\mu^{2}},\frac{x}{\mu}\right)\]
is a solution for the IVP (\ref{givp1}) with initial data
$w^{\mu}_{0}=\mu^{-1/2}w(\mu^{-1}x)$.
 From Plancherel's theorem and a simple computation we see that
$$
      \|I\partial_{x}w^{\mu}_{0}\|_2\lesssim
      \frac{N^{1-s}}{\mu^{s}}\|w_{0}\|_{H^{s}}.
$$
while
$$ \| I w^{\mu}_0 \|_2 \leq \| w^{\mu}_0\|_2 = \| w_0 \|_2 <
\sqrt{2\pi}.$$
We now choose $\mu := N^{\frac{1-s}{s}}$.  From the previous we see
that
$\|Iw^{\mu}_{0}\|_{H^1}\lesssim 1$,
so from Sobolev embedding (or Gagliardo-Nirenberg) we obtain
$$ E(Iw^\mu_0) \leq C_1$$
for some constant $C_1 > 0$.

Now suppose inductively that we have a time $T$ such that
$$ E(Iw^\mu(T)) \leq C_1 + C_2 N^{-1+} T$$
where $C_2 > 0$ is a constant depending on $C_1$ to be chosen later.
If $T \ll N^{1-}$, we then have $E(Iw^\mu(T)) \leq 2C_1$, which implies
from Lemma \ref{gagliardo} that
$$ \| Iw^\mu(T) \|_{H^1} \leq C_3$$
where $C_3$ depends on $C_1$.  By Theorem \ref{lwp} we thus have
$$ \| Iw^\mu \|_{X^{1,1/2+}([T,T+\delta] \times \R)} \leq C_4$$
where $C_4$, $\delta$ depend on $C_3$.

In the next four sections we shall prove the key estimate

\begin{lemma}\label{pinc}  For any Schwartz function $w$, we have
\begin{equation}\label{pinc-est}
|\int_T^{T+\delta} \Lambda_n(M_n; w(t))\ dt|
\lesssim N^{-1+} \| Iw\|_{X^{1,1/2+}([T,T+\delta]\times \dbR)}^n
\end{equation}
for $n=4,6,8$, where $M_4$, $M_6$, $M_8$ are defined in Proposition
\ref{energy-increment}.
\end{lemma}

Assuming this estimate for the moment, we see from the previous and
Proposition \ref{energy-increment} that
$$ E(Iw^\mu(T+\delta)) \leq E(Iw^\mu(T)) + C_5 N^{-1+}$$
where $C_5$ depends on $\delta$ and $C_4$.  This allows us to close
the induction hypothesis by setting $C_2 := C_5$.  As a consequence
we have thus shown that\footnote{Strictly speaking, we have only
shown this for $T$ being an integer multiple of $\delta$, however
this can be easily remedied, e.g. by using the fact that the
$X^{1,1/2+}$ norm controls the $L^\infty_t H^1_x$ norm on
$[T,T+\delta] \times \R$.}
$$ \| I w^\mu(T) \|_{H^1} \lesssim 1$$
for all $T \ll N^{1-}$.  From the definition of $I$ this implies that
$$ \| w^\mu(T) \|_{H^s} \lesssim C_N$$
for all $T \ll N^{1-}$.  Undoing the scaling, this implies that
$$ \| w(T) \|_{H^s} \lesssim C_{N,\mu}$$
for all $T \ll N^{1-}/\mu^2$.  However, if $s > 2/3$, then
$N^{1-}/\mu^2 = N^{\frac{3s-2-}{s}}$ goes to infinity as $N \to \infty$,
and Proposition \ref{global} follows.
\hfill$\square$

\begin{remark}
An examination of the above argument shows also that the $H^s$ norm
of $w$ (and of $u$) grows at most polynomially in time, however the
order of this growth obtained by this argument goes to infinity as $s
\to 2/3$.
\end{remark}

\section{Proof of Lemma \ref{pinc}: preliminaries}\label{prelim-sec}

To prove Lemma \ref{pinc} we shall treat the cases $n=4$, $n=6$,
$n=8$ separately.  The idea will be first to obtain some good
estimates on $M_n$ in terms of $m(\xi_i)$ and $\langle \xi_i
\rangle$, and then to bound the resulting multilinear expression
using standard tools such as the Strichartz estimates
\eqref{strichartz-6}, \eqref{strichartz-infty}, the trivial estimate
\begin{equation}\label{triv}
\| u \|_{L^2_t L^2_x} \lesssim \| u \|_{0,0},
\end{equation}
and H\"older's inequality.  In addition to the above linear
estimates, we shall also take advantage of the following bilinear
improvement to Strichartz' estimate in the case of differing
frequencies (cf. \cite{borg:refinements})

\begin{lemma}\label{improved-strichartz}
For any Schwartz functions $u, v$ with Fourier support in $|\xi| \sim
R$, $|\xi| \ll R$ respectively, we have that
$$ \| u v \|_{L^2_t L^2_x} = \| u \bar v \|_{L^2_{t} L^2_x}
\lesssim R^{-1/2} \|u\|_{0,1/2+} \|v\|_{0,1/2+}.$$

\end{lemma}

\begin{proof}
This is an improved Strichartz estimate of the type considered in
\cite{borg:refinements}.

It is enough to show that  if $u$ and $v$ are solutions of the
free Schr\"odinger equation, that is $u=e^{it\partial_{x}^{2}}\phi$
and
   $v=e^{it\partial_{x}^{2}}\psi$, then
\begin{equation}
\|D_{x}^{1/2}(uv)\|_{L^{2}}\lesssim\|\phi\|_{L^{2}}\|\psi\|_{L^{2}}.
\label{stl},\end{equation}
where $D_{x}$ is the operator such that $\widehat{D_{x}f}(\xi)
=\langle\xi\rangle\hat{f}(\xi)$.
If we use duality and the change of variable $ \xi_{1}+\xi_{2}=s$ and
$|\xi_{1}|^{2}+|\xi_2|^{2}=r$, the left hand side of \eqref{stl}
becomes
\begin{eqnarray*}
& &\sup_{\|F\|_{L^{2}}\leq 1}\int R^{1/2} F(\xi_{1}+\xi_{2},
|\xi_{1}|^{2}+|\xi_{2}|^{2})\hat{\phi}(\xi_{1})
\hat{\psi}(\xi_{2})d\xi_{1}d\xi_{2}\\
&\lesssim&\int R^{1/2} F(s,r)\frac{H(s,r)}{R}ds dr,
\end{eqnarray*}
where $H(s,r)$ denotes the product of $\hat{\phi}$ and 
$\hat{\psi}$ in the new
variables. Notice that the change of variables introduced above has a
Jacobian of size $R$. Now if we use Cauchy-Schwarz and we change the
variables
back to $\xi_{1}$ and $\xi_{2}$, we obtain (\ref{stl}).
\end{proof}

In one of our sub-cases, we shall also take advantage of a trick
(originally due to Bourgain \cite{borg:xsb}) of splitting the symbol $|\xi_1|^2
- \ldots + |\xi_n|^2$ as a sum of $\tau_j \mp |\xi_j|^2$.

Our estimates are not best possible, and it is likely that one can
improve the $N^{-1+}$ gain in our estimates, probably to
$N^{-3/2+}$.  However this will fall short of the $N^{-2+}$ gain
needed to push the global well-posedness down to match the local
well-posedness theory at $s > 1/2$.  However one can  recover
this by adding higher order correction terms to the energy
$E_N(w(t))$, as in \cite{ckstt:2}.
If one does this, one will end up estimating $\Lambda_6$ and
$\Lambda_8$ expressions rather than $\Lambda_4$.  This will be
beneficial because such expressions will have fewer derivatives in
their symbol and can therefore enjoy better decay in $N$. The details
of this argument will appear in a later paper.

We set out some notation.  Let $n = 4$, $6$, or $8$, and let $\xi_1,
\ldots, \xi_n$ be frequencies such that $\xi_1 + \ldots + \xi_n =
0$.  Define $N_i := |\xi_i|$, and $N_{ij} := |\xi_{ij}|$.  We adopt
the notation that
$$ 1 \leq soprano, alto, tenor, baritone \leq n$$
are the distinct indices such that
$$ N_{soprano} \geq N_{alto} \geq N_{tenor} \geq N_{baritone}$$
are the highest, second highest, third highest, and fourth highest
values of the frequencies $N_1, \ldots, N_n$ respectively (if there
is a tie in frequencies, we break the tie arbitrarily).

Since $\xi_1 + \ldots + \xi_n = 0$, we must have $ N_{soprano} \sim
N_{alto}$.  Also, from Proposition \ref{energy-increment} we see that
$M_n$ vanishes unless $N_{soprano} \gtrsim N$.

\section{Proof of Lemma \ref{pinc} when $n=4$}\label{proof-4}
We now estimate the $\Lambda_4$ expression.  We begin by estimating
the multiplier $M_4$.

\begin{lemma}\label{m4-est}  Let $\xi_1,\xi_2,\xi_3,\xi_4$ be such
that $\xi_{1234} = 0$.
\begin{itemize}
\item  If $N_{tenor} \sim N_{soprano}$, then
\begin{equation}\label{high-tenor}
|M_4(\xi_1,\xi_2,\xi_3,\xi_4)| \lesssim
N^{-1} (N/N_{soprano})^{1/10} \langle \xi_{12} \xi_{14} \rangle^{1/2}
\prod_{j=1}^4 \langle \xi_j \rangle m(\xi_j).
\end{equation}
\item  If $N_{tenor} \ll N_{soprano}$, then
\begin{equation}\label{low-tenor}
|M_4(\xi_1,\xi_2,\xi_3,\xi_4)| \lesssim
N^{-1} (N/N_{soprano})^{1/10} N_{soprano} \prod_{j=1}^4 \langle \xi_j
\rangle m(\xi_j).
\end{equation}
\end{itemize}
\end{lemma}

\begin{proof}  Fix $\xi_1, \ldots, \xi_4$.
If $N_{soprano} \ll N$ then $M_4$ vanishes by the second part of
Proposition \ref{energy-increment}, so we will assume that
$N_{soprano} \gtrsim N$.

We split $M_4 = C_1 M'_4 + C_2 M''_4$, where
$$
M'_4 := m_1 m_2 m_3 m_4 \xi_{12} \xi_{13} \xi_{14}
$$
and
$$
M''_4 := m_1^2 \xi_1^2 \xi_3 + m_2^2 \xi_2^2 \xi_4 + m_3^2 \xi_3^2
\xi_1 + m_4^2 \xi_4^2 \xi_2.
$$
In the $N_{soprano} \gtrsim N$ case we will not need to exploit
cancellation between $M'_4$ and $M''_4$ (although such cancellation
certainly exists), and shall estimate them separately.

Let us first prove \eqref{high-tenor}.  We begin with estimating
$M'_4$.  We have
\begin{align*}
|M'_4|
&= N_{12} N_{13} N_{14} m(N_1) m(N_2) m(N_3) m(N_4)\\
&\lesssim N_{12} N_{14} N_{soprano} m(N_{soprano})^3\\
&\lesssim \langle N_{12} N_{14} \rangle^{1/2} N_{soprano}^2
m(N_{soprano})^3\\
&\lesssim \langle N_{12} N_{14} \rangle^{1/2} \frac{1}{N}
(N/N_{soprano})^{1/10} \langle N_{soprano} \rangle^3
m(N_{soprano})^3\\
&\sim N^{-1} (N/N_{soprano})^{1/10} \langle N_{12} N_{14}
\rangle^{1/2}  \langle N_{soprano}\rangle m(N_{soprano}) \langle
N_{alto} \rangle m(N_{alto}) \langle N_{tenor} \rangle m(N_{tenor})\\
&\lesssim N^{-1} (N/N_{soprano})^{1/10} \langle N_{12} N_{14}
\rangle^{1/2} \prod_{j=1}^4 \langle N_j \rangle m(N_j)
\end{align*}
as desired.

It remains to estimate $M''_4$.  We divide into two cases:
$N_{baritone} \sim N_{soprano}$ and $N_{baritone} \ll N_{soprano}$.

{\bf Case 1:  $N_{baritone} \sim N_{soprano}$.}

In this case all the frequencies are comparable to each other.  By
symmetry we may assume that $N_{12} \leq N_{14}$, in which case it
suffices to show
$$ |M''_4| \lesssim N^{-1} (N/N_1)^{1/10} N_{12} N_1^4 m(N_1)^4.$$
We can rewrite $M''_4=f(0) - f(h)$, where
$$ f(h) := m(\xi_1-h)^2 (\xi_1-h)^2 (\xi_3+h) + m(\xi_3+h)^2
(\xi_3+h)^2 (\xi_1-h)$$
and $h := \xi_1 + \xi_2$.  A routine calculation shows that
$$ |f'(x)| \lesssim m(N_1)^2 N_1^2$$
for all $x = O(N_1)$, so by the mean value theorem and the assumption
$N_1 \gtrsim N$ we have
$$ |M''_4| = |f(0)-f(h)| \lesssim N_{12} m(N_1)^2 N_1^2 \lesssim
N^{-1} (N/N_1)^{1/10} N_{12} N_1^4 m(N_1)^4$$
as desired (in fact we gain an additional power of $N$).

{\bf Case 2:  $N_{baritone} \ll N_{soprano}$.}

By symmetry we may assume that $baritone = 4$, thus $N_1 \sim N_2
\sim N_3 \gg N_4$.  In this case $N_{14} \sim N_1$, $N_{12} = N_{34}
\sim N_1$ and $\langle N_4 \rangle m(N_4) \gtrsim 1$, so it suffices
to show
$$ |M''_4| \lesssim N^{-1} (N/N_1)^{1/10} N_1^4 m(N_1)^3.$$
But we may crudely estimate the left-hand side by
$$ |M''_4|\lesssim  m(N_1)^2 N_1^3 + m(N_4)^2 N_4^2 N_1 \lesssim
m(N_1)^2 N_1^3$$
which suffices since $N_1 \gtrsim N$.  This proves \eqref{high-tenor}.

Now we show \eqref{low-tenor}. Observe that
$$ N_{12} N_{13} N_{14} \lesssim N_{soprano}^2 N_{tenor}$$
and hence
\begin{align*}
|M'_4|
&\lesssim N_{soprano}^2 N_{tenor} m(N_{soprano}) m(N_{alto})
m(N_{tenor}) m(N_{baritone})\\
&\lesssim \prod_{j=1}^4 \langle N_j \rangle m(N_j)\\
&\lesssim N^{-1+} N_{soprano}^{1-} \prod_{j=1}^4 \langle N_j \rangle
m(N_j).
\end{align*}
Thus it only remains to estimate $M''_4$.  Since $\langle
N_{baritone} \rangle m(N_{baritone})$ and \\
$ m(N_{tenor}) N^{-1}
(N/N_{soprano})^{1/10} N_{soprano}$ are both $\gtrsim 1$, it suffices
to show
$$ |M''_4| \lesssim m(N_{soprano})^2 N_{soprano}^2 N_{tenor}.$$

By symmetry we may reduce to one of two cases.

{\bf Case 1:  $N_3 = N_{tenor}$ and $N_4 = N_{baritone}$.}

We  crudely estimate
$$ |M''_4(\xi_1,\xi_2,\xi_3,\xi_4)| \lesssim
m(N_1)^2 N_1^2 N_3 + m(N_2)^2 N_2^2 N_4 + m(N_3)^2 N_3^2 N_1
+ m(N_4)^2 N_4^2 N_2 \lesssim m(N_1)^2 N_1^2 N_3$$
as desired.

{\bf Case 2:  $N_2 = N_{tenor}$ and $N_4 = N_{baritone}$.}

In this case we estimate
\begin{align*}
|M''_4| &\lesssim
|m_1^{2} \xi_1^2 \xi_3 + m_3^2 \xi_3^2 \xi_1| +
m(N_2)^2 N_2^2 N_4 + m(N_4)^2 N_4^2 N_2\\
&=
N_1 N_3 |m(\xi_1)^2 \xi_1 - m(\xi_1 + \xi_2 + \xi_4)^2 (\xi_1 + \xi_2
+ \xi_4)|
+ O(m(N_1)^2 N_1^2 N_2).
\end{align*}
The function $m(\xi_1 + h)^2 (\xi_1+h)$ has a derivative of
$O(m(N_1)^2)$
whenever $|h| \ll N_{1}$, thus by the mean value theorem we thus have
$$
|M'_4(\xi_1,\xi_2,\xi_3,\xi_4)| \lesssim
N_1 N_3 N_{24} m(N_1)^2
+ O(m(N_1)^2 N_1^2 N_2)
\lesssim m(N_1)^2 N_1^2 N_2$$
as desired.
\end{proof}

We now prove \eqref{pinc-est} in the $n=4$ case. It suffices to show
that
$$
\int_T^{T+\delta} \Lambda_4(M_4; w_1(t), \overline{w_2(t)}, w_3(t),
\overline{w_4(t)})\ dt \lesssim N^{-1+}
\prod_{j=1}^4 \| I w_j\|_{1,1/2+}$$
for all Schwartz functions $w_1, \ldots, w_4$ on $\R \times \R$.
Since $M_4$ vanishes for $N_{soprano} \ll N$, it  suffices by
dyadic decomposition to show that
$$
\int_T^{T+\delta} \Lambda_4(M_4 \chi_{N_{soprano} \sim 2^k}; w_1(t),
\overline{w_2(t)}, w_3(t), \overline{w_4(t)})\ dt \lesssim N^{-1+}
2^{(0+)k} (N/2^k)^{1/10} \prod_{j=1}^4 \| I w_j\|_{1,1/2+}$$
for all integers $k$ for which $2^k \gtrsim N$. (The exact choice of
the cutoff $\chi_{N_{soprano} \sim 2^k}$ is not important as we shall
soon be taking absolute values everywhere anyway).

Fix $k$. Without loss of generality we may assume that the Fourier
transforms $\tilde w_j$ are real and non-negative.  We divide into
the $\Lambda_4$ integral into the regions $N_{tenor} \sim
N_{soprano}$ and $N_{tenor} \ll N_{soprano}$.

{\bf Case 1. $N_{tenor} \sim N_{soprano}$.}

We first perform some manipulations to eliminate the cutoff
$\chi_{[T,T+\delta]}(t)$.  Write $\chi_{[T,T+\delta]}(t) = a(t) +
b(t)$, where $a(t)$ is $\chi_{[T,T+\delta]}(t)$ convolved with a
smooth approximation to the identity of width $2^{-100k}$, and $b(t)
= \chi_{[T,T+\delta]}(t) - a(t)$.

Let us first consider the contribution of $b(t)$.  We crudely
estimate $M_4 = O(2^{10k})$ and estimate this contribution by
$$ 2^{10k} \int\int |b(t)| |w_1(t,x)| |w_2(t,x)| |w_3(t,x)|
|w_4(t,x)|\ dx dt.$$
By H\"older, three applications of \eqref{strichartz-6}, one
application of \eqref{strichartz-2}, and four applications of
\eqref{i-smoothing} we can bound this by
$$ 2^{10k} \|b\|_2 \prod_{j=1}^4 \| Iw_j \|_{1,1/2+}.$$
Since $\|b\|_2 \lesssim 2^{-50k}$, the claim then follows.

Now consider the contribution of $a(t)$.  We use

\begin{lemma}\label{a-bound}
We have
$$ \| a(t) w_1 \|_{1,1/2+} \lesssim 2^{(0+)k} \| w_1 \|_{1,1/2+}.$$
\end{lemma}

\begin{proof}
By applying Plancherel, restricting to a single frequency $\xi$, and
then undoing Plancherel, we see that it suffices to show that
$$ \| a(t) f \|_{H^{1/2+}_t} \lesssim 2^{(0+)k} \| f \|_{H^{1/2+}_t}$$
for all functions $f$.  But this follows from the routine calculation
$$ \| a(t) \|_{H^{1/2+}_t} \lesssim 2^{(0+)k}$$
and the fact that $H^{1/2+}_t$ is closed under multiplication.
\end{proof}

It therefore suffices to show
$$
|\int \Lambda_4(M_4 \chi_{N_{tenor} \sim N_{soprano} \sim 2^k};
w_1(t), \overline{w_2(t)}, w_3(t), \overline{w_4(t)})\ dt| \lesssim
N^{-1} (N/2^k)^{1/10}
\prod_{j=1}^4 \| I w_j\|_{1,1/2+}.$$
Without loss of generality we may assume that the Fourier transforms
$\tilde w_j$ are real and non-negative.  By Plancherel and
\eqref{high-tenor} we estimate the left-hand side by
$$
N^{-1} (N/2^k)^{1/10} |\int_*  \langle \xi_{12} \xi_{14}
\rangle^{1/2}
\widetilde{ID_x w_1}(\tau_1,\xi_1)
\widetilde{\overline{ID_x w_2}}(\tau_2,\xi_2)
\widetilde{ID_x w_3}(\tau_3,\xi_3)
\widetilde{\overline{ID_x w_4}}(\tau_4,\xi_4)|.$$
 From the identity (cf. Bourgain \cite{borg:xsb} and Kenig-Ponce-Vega
\cite{kpv:kdv})
\begin{align*}
\sum_{j=1}^4 (\tau_j - (-1)^{j-1} \xi^2_j) &=
-\xi_1^2 + \xi_2^2 - \xi_3^2 + \xi_4^2\\
&= \xi_{12} \xi_{2-1} + \xi_{34} \xi_{4-3}\\
&= \xi_{12} (\xi_{2-1} - \xi_{4-3})\\
&= - 2 \xi_{12} \xi_{14}
\end{align*}
we see that
$$ \langle \xi_{12} \xi_{14} \rangle \lesssim \langle \tau_j -
(-1)^{j-1} \xi^2_j \rangle$$
for some $j=1,2,3,4$.  We shall assume $j=1$; the argument for other
values of $j$ is similar.  We can then use duality and Plancherel to
estimate the previous by
$$
N^{-1} (N/2^k)^{1/10} \| I w_1 \|_{1,1/2+}
\| \overline{ID_x w_2} ID_x w_3 \overline{ID_x w_4} \|_{L^2_t
L^2_x}.$$
But this is acceptable by H\"older and three applications of
\eqref{strichartz-6}.

{\bf Case 2. $N_{tenor} \ll N_{soprano}$.}

We shall assume that $soprano = 1$ and $alto = 2$; the reader may
verify that the other cases follow by the same argument.  We may then
restrict $w_1$, $w_2$ to have Fourier support in $|\xi| \sim 2^k$ and
$w_3$, $w_4$ to have Fourier support in the region $|\xi| \ll 2^k$.

By \eqref{low-tenor} we have
$$ |M_4| \lesssim N^{-1} (N/2^k)^{1/10}  2^k \prod_{j=1}^4 \langle
\xi_i \rangle m(\xi_i).$$
The claim then follows from H\"older and two applications of
Proposition \ref{improved-strichartz}.

\section{Proof of Lemma \ref{pinc} when $n=6$}\label{proof-6}

We begin with the analogue of Lemma \ref{m4-est}.

\begin{lemma}\label{m6-est}
Let $\xi_1, \ldots, \xi_6$ be such that $\xi_{123456}=0$.
\begin{itemize}
\item  If $N_{tenor} \sim N_{soprano}$, then
\begin{equation}\label{high-six}
|M_6(\xi_1,\xi_2,\xi_3,\xi_4,\xi_5,\xi_6)| \lesssim
N^{-1}
\langle \xi_{soprano} \rangle m(\xi_{soprano})
\langle \xi_{alto} \rangle m(\xi_{alto})
\langle \xi_{tenor} \rangle m(\xi_{tenor}).
\end{equation}
\item  If $N_{tenor} \ll N_{soprano}$, then
\begin{equation}\label{low-six}
|M_6(\xi_1,\xi_2,\xi_3,\xi_4,\xi_5,\xi_6)| \lesssim
N^{-1+} \langle N_{soprano} \rangle^{1-}
\langle \xi_{soprano} \rangle m(\xi_{soprano})
\langle \xi_{alto} \rangle m(\xi_{alto}).
\end{equation}
\end{itemize}
\end{lemma}

One can improve these estimates by exploiting further cancellation in
the expression $M_6$, but we shall not need to do so because of the
good smoothing properties of our equation \eqref{givp1}.

\begin{proof}
Since $\xi_{123456}=0$, we have $N_{alto} \sim N_{soprano}$.  We may
also assume that $N_{soprano} \gtrsim N$ since $M_6$ vanishes
otherwise.

We have the very crude estimate
$$ |M_6| \lesssim N_{soprano}^2.$$

If $N_{tenor} \sim N_{soprano}$, we then have
$$ |M_6| \lesssim N_{soprano}^2 \lesssim N^{-1} m(N_{soprano})
N_{soprano} m(N_{alto}) N_{alto} m(N_{tenor}) N_{tenor}$$
(using the hypothesis $s>2/3$), and \eqref{high-six} follows.

Now suppose that $N_{tenor} \ll N_{soprano}$.  Then
$$ |M_6| \lesssim N_{soprano}^2 \lesssim N^{-1+} \langle N_{soprano}
\rangle^{1-}  m(N_{soprano}) N_{soprano} m(N_{alto}) N_{alto}$$
(since $s > 1/2$), and \eqref{low-six} follows.
\end{proof}

We now prove \eqref{pinc-est} for $n=6$.  As in the previous section
it suffices to show
$$
\int_T^{T+\delta} \Lambda_6(M_6; w_1(t), \overline{w_2(t)}, w_3(t),
\overline{w_4(t)}, w_5(t), \overline{w_6(t)})\ dx dt \lesssim N^{-1+}
\prod_{j=1}^6 \| I w_j\|_{1,1/2-}$$
for all Schwartz functions $w_1, \ldots, w_6$ on $\R \times \R$.
Without loss of generality we may assume that the Fourier transforms
$\tilde w_i$ of $w_i$ are real and non-negative.

We again divide into the cases $N_{tenor} \sim N_{soprano}$ and
$N_{tenor} \ll N_{soprano}$.

{\bf Case 1. $N_{tenor} \sim N_{soprano}$.}

By \eqref{high-six} and symmetry it suffices to show
$$
\int_T^{T+\delta} \int \prod_{j=1}^3 |D_x I w_j|  \prod_{j=4}^6
|w_j|\ dx dt \lesssim
\prod_{j=1}^6 \| I w_j\|_{1,1/2+}.$$
But this follows from H\"older, six applications of
\eqref{strichartz-6} first, and three applications of
\eqref{i-smoothing} after.

{\bf Case 2. $N_{tenor} \ll N_{soprano}$.}

We shall assume that $soprano = 1$ and $alto = 2$; the reader may
verify that the other cases follow by the same argument.

First suppose that $N_{soprano} \sim 2^k$ for some integer $k$.  Then
$w_1, w_2$ have Fourier support on $|\xi| \sim 2^k$, while $w_3, w_4,
w_5, w_6$ have Fourier support on $|\xi| \ll 2^k$.

We apply \eqref{low-six}, and bound the contribution of this case by
$$
N^{-1+} 2^{(1-)k}
\int_T^{T+\delta} \int \prod_{j=1}^2 |D_x I w_j|  \prod_{j=3}^6
|w_j|\ dx dt,$$
which we bound using H\"older by
$$
N^{-1+} 2^{(1-)k}
\| (D_x I w_1) w_3 \|_{L^2_t L^2_x}
\| (D_x I w_2) w_4 \|_{L^2_t L^2_x}
\| w_5 \|_{L^\infty_t L^\infty_x}
\| w_6 \|_{L^\infty_t L^\infty_x}.$$
By Lemma \ref{improved-strichartz}, \eqref{strichartz-infty}, and
\eqref{i-smoothing} we can bound this by
$$ N^{-1+} 2^{(0-)k} \prod_{j=1}^6 \| I w_j\|_{1,1/2+}$$
The claim then follows by summing in $k$.

\section{Proof of Lemma \ref{pinc} when $n=8$}\label{proof-8}

We begin with the analogue of Lemma \eqref{m4-est}.

\begin{lemma}\label{m8-est}
For any $\xi_1, \ldots, \xi_6$ with $\xi_{123456}=0$, we have
\begin{equation}\label{eights}
|M_8(\xi_1,\xi_2,\xi_3,\xi_4,\xi_5,\xi_6,\xi_7,\xi_8)| \lesssim
N^{-1}
\langle \xi_{soprano} \rangle m(\xi_{soprano})
\langle \xi_{alto} \rangle m(\xi_{alto}).
\end{equation}
\end{lemma}

\begin{proof}
As usual we may assume that $N_{soprano} \sim N_{alto} \gtrsim N$.
We crudely estimate
$$ |M_8| \lesssim N_{soprano} \lesssim N^{-1} N_{soprano}
m(N_{soprano}) N_{alto} m(N_{alto})$$
and the claim follows.
\end{proof}

To prove \eqref{pinc-est} for $n=8$ it suffices to show
$$
\int_T^{T+\delta} \Lambda_8(M_8; w_1(t), \ldots, \overline{w_8(t)})\
dt \lesssim N^{-1}
\prod_{j=1}^8 \| I w_j\|_{1,1/2+}$$
for all Schwartz functions $w_1, \ldots, w_8$ on $\R \times \R$.
Without loss of generality we may assume that the Fourier transforms
$\tilde w_i$ of $w_i$ are real and non-negative.  By Lemma
\ref{m8-est} and symmetry it thus suffices to show
$$
\int_T^{T+\delta}\int |D_x I w_1| |D_x I w_2| \prod_{j=3}^8 |w_j|\ dx
dt \lesssim N^{-1}
\prod_{j=1}^8 \| I w_j\|_{1,1/2+}.$$
But this follows from H\"older, six applications of
\eqref{strichartz-6}, and two applications of
\eqref{strichartz-infty} and \eqref{i-smoothing}.

\begin{remark}
As it was shown in Section 3, the gauge transform  in
Definition \ref{gauge-def} introduces a quintic term  in the initial
value problem (\ref{givp1}). Then one can ask if  the same arguments
we proposed above can be used in order to study the global well-posedness
of the quintic nonlinear Schr\"odinger initial value problem
\begin{equation}
\left\{ \begin{array}{l}
i\partial_tv+\partial_{x}^{2}v +\lambda|v|^{4}v=0,\\
    v(x,0) = v_{0}(x),\hspace{1.5cm}x \in \dbR, \, t \in \dbR,
\end{array}\right.
\label{quintic}\end{equation}
where $\lambda \in \dbR$. In this case we define the energy
$$H(f) := \int |\partial_x f(x)|^2\ dx - \frac{\lambda}{6}\int |f|^{2}\ dx.$$
By Plancherel, we may write $H(f)$ using the $\Lambda$ notation as
$$ H(f) = \Lambda_2(\xi_1\xi_2;f) - \frac{\lambda}{6}\Lambda_6(1;f).$$
As in Lemma \ref{conservation} one can prove that the energy $H(v(t))$
of the solution $v$ for (\ref{quintic}) is constant. Now let's define
the new energy
$$H_{N}(v)=H(Iv)=\Lambda_2(\xi_1\xi_2m_{1}m_{2};v)
-\frac{\lambda}{6}\Lambda_6(\prod_{i=1}^{6}m_{i};v).$$
just like we did in Secition 4. Then by the analogue of (\ref{diff})
$\partial_{t}H_{N}(v(t))$ will involve terms of type $\Lambda_{2},
\Lambda_{6}$ and $\Lambda_{10}$. Using the same ideas
presented in the proof of  Lemma \ref{pinc} we can estimate in the
appropriate way also the term involving $\Lambda_{10}$.
If in (\ref{quintic}) we assume that $\lambda<0$ (defocusing) or that
the $L^{2}$ norm of the initial data is small (so that the
Gagliardo-Niremberg inequality can be applied) then the  energy
$H(v)(t)$ stays positive for all times and
global well-posedness in $H^{s}$ for $s>2/3$ will follow.
We will present the details of the
proof in a future paper. It has to be said here that global results
for ``small data'' are already available for (\ref{quintic})
through more standard arguments \cite{cw}.
\end{remark}

\end{document}